\documentclass[a4paper,11pt]{article}
\usepackage{fontspec}
\usepackage[english]{babel}
\usepackage{amsmath,amssymb,amsthm,mathtools,mathrsfs}
\usepackage[margin=25mm,headheight=15pt]{geometry}
\usepackage{microtype}
\usepackage{titling}
\usepackage{booktabs,longtable,array,calc}
\usepackage{enumitem,xcolor,xurl}
\usepackage{tikz-cd}
\usepackage[numbers,square,sort&compress]{natbib}
\usepackage[unicode]{hyperref}
\usepackage{bookmark}
\hypersetup{
 colorlinks=true,linkcolor=blue!45!black,citecolor=blue!45!black,
 urlcolor=blue!55!black,
 pdftitle={Higher Noncommutative Motives and Higher K-theory},
 pdfsubject={Additive motives from coefficient modules in infinity categories},
 pdfcreator={XeLaTeX}}
\setlist{itemsep=0.2em,topsep=0.4em,parsep=0.1em}

\allowdisplaybreaks[1]
\numberwithin{equation}{section}
\theoremstyle{plain}
\newtheorem{theorem}{Theorem}[section]
\newtheorem{lemma}[theorem]{Lemma}
\newtheorem{proposition}[theorem]{Proposition}
\newtheorem{corollary}[theorem]{Corollary}

\theoremstyle{definition}
\newtheorem{definition}[theorem]{Definition}

\theoremstyle{remark}
\newtheorem{remark}[theorem]{Remark}
\newtheorem{example}[theorem]{Example}

\title{Higher Noncommutative Motives and Higher K-theory}
\author{Hisato Matsukawa\thanks{Department of Mathematics, Faculty of Science, Hokkaido University\\Kita 10, Nishi 8, Kita-Ku, Sapporo, Hokkaido, 060-0810, Japan\\Email:matsukawa.hisato.f4@elms.hokudai.ac.jp,hmnr0211@gmail.com}}
\date{}
\makeatletter
\ifx\@author\@empty
  \preauthor{}\postauthor{}
\fi
\ifx\@date\@empty
  \predate{}\postdate{}
\fi
\makeatother
\begin{document}
\maketitle
\begin{abstract}
We construct higher additive noncommutative motives of stable
$(\infty,n)$-categories, extending the additive motives of
Blumberg, Gepner, and Tabuada and the secondary theory of
Mazel-Gee and Stern.
We formulate
a notion of exact sequences and establish the universal
property of the resulting additive motives.
We give two equivalent constructions: a localization of
enriched presheaves and a weighted Waldhausen
$S$-construction on complete exact flags.
We then compute their Hom objects, obtaining ordinary
K-theory at the first level and lower additive motives
recursively at higher levels.
Our results recover the classical and secondary theories,
provide coherent multiplicative structures, and yield
an endofunctor on the category of Stefanich rings.
In short, higher motives are higher K-theory.

\end{abstract}
\newcommand{\Sp}{\mathrm{Sp}}
\newcommand{\Cat}{\operatorname{Cat}}
\newcommand{\PrL}{\mathrm{Pr}^{L,\omega}}
\newcommand{\nPr}[1]{{#1}\mathrm{Pr}_{\mathrm{st}}}
\newcommand{\M}{\mathsf{Mot}}
\newcommand{\PrHom}{\mathcal{H}\!om}
\newcommand{\Map}{\operatorname{Map}}
\newcommand{\Fun}{\operatorname{Fun}}
\newcommand{\Mod}{\operatorname{Mod}}
\newcommand{\Perf}{\operatorname{Perf}}
\newcommand{\EqRow}{\operatorname{ExSeq}}
\newcommand{\cofib}{\operatorname{cofib}}
\newcommand{\id}{\mathrm{id}}
\newcommand{\Ind}{\operatorname{Ind}}
\newcommand{\Add}{\operatorname{Add}}
\newcommand{\core}{\simeq}
\newcommand{\ex}{\mathrm{ex}}
\newcommand{\op}{\mathrm{op}}
\newcommand{\Pmod}{\operatorname{PSh}}
\newcommand{\one}{\mathbf 1}

\section{Introduction}

The universal characterization of algebraic K-theory by
Blumberg, Gepner, and Tabuada provides a systematic approach
to additive invariants of stable infinity-categories
\citep{BGT}. Mazel-Gee and Stern extended this perspective
to stable $(\infty,2)$-categories, obtaining category-valued
additive motives and secondary algebraic K-theory \citep{MGS}.
In this paper, we develop a theory of higher additive motives
that extends these constructions to stable $(\infty,n)$-categories.

Our framework is based on the symmetric monoidal tower
\[
 \nPr0=\Sp,\qquad
 \nPr{(r+1)}=\Mod_{\nPr r}(\PrL),
\]
whose precise definition is given in
Section~\ref{sec:conventions}.
The relevant exact sequences generalize cofiber sequences
and the localization sequences of Mazel-Gee and Stern,
incorporating adjunctions and exactness in the Hom objects.

Our first main result constructs universal additive
invariants at every finite height.

\begin{theorem}[Additive motives]
\label{thm:main}
For every $n\geq1$ and $A\in\nPr{(n+1)}$,
there exist $\M_n(A)\in\nPr n$ and a canonical
additive invariant $U_A:A\to \M_n(A)$ such that
\begin{equation}\label{eq:main-up}
 \Fun^L_{\nPr{(n-1)}}(\M_n(A),E)
 \simeq \Add_{n,\mathrm{filt}}(A,E)
\end{equation}
naturally in $E\in\nPr n$.
\end{theorem}

The right-hand side, $\Add_{n,\mathrm{filt}}(A,E)$, is the category of additive invariants valued in $E$.

We give two constructions of these motives.
The first is obtained by localizing enriched presheaves
(Theorem~\ref{thm:presheaf-universal}),
generalizing the constructions of additive motives in
\cite{BGT,MGS}. The second is a higher analogue of
Waldhausen's $S$-construction
(Theorem~\ref{thm:weighted-s-universal}).
A key ingredient is the existence of compact modules
classifying exact sequences.
We establish the universal properties of both constructions,
thereby proving their equivalence.

Having constructed the motives, we compute their Hom objects.
This reveals a recursive relationship between additive motives
at successive heights.

\begin{theorem}[Mapping objects;
Theorem~\ref{thm:mapping-comparison}]
\label{thm:intro-mapping}
Let $A\in\nPr{(n+1)}$, let $c\in A$ be compact,
and let $a\in A$ be arbitrary. Then
\[
 \operatorname{map}_{\M_1(A)}(U_Ac,U_Aa)
 \simeq K\bigl(\PrHom_A(c,a)^\omega\bigr)
 \qquad(n=1),
\]
and
\[
 \PrHom_{\M_n(A)}(U_Ac,U_Aa)
 \simeq \M_{n-1}\bigl(\PrHom_A(c,a)\bigr)
 \qquad(n\geq2),
\]
where $K$ denotes connective algebraic K-theory.
These equivalences are compatible with composition,
units, and external products.
\end{theorem}

In other words, the higher additive motives themselves
are the values of higher K-theory.

These results recover the additive motives of \cite{BGT,HSS}
over rigid stable bases and the category-valued
secondary motives of Mazel-Gee and Stern.

Finally, we study the interaction of higher K-theory
with Stefanich rings. 
For $A\in\nPr n$, higher K-theory is defined
by
\[
K^{[n]}(A)=
\begin{cases}
\operatorname{map}_{\M_1(\nPr1)}
\bigl(U_{\nPr1}(\Sp),U_{\nPr1}(A)\bigr),
&n=1,\\
\PrHom_{\M_n(\nPr n)}
\bigl(U_{\nPr n}(\nPr{(n-1)}),U_{\nPr n}(A)\bigr),
&n\geq2.
\end{cases}
\]
The multiplicative compatibility of our
mapping formula yields the following result.

\begin{theorem}[K-theory of Stefanich rings;
Theorem~\ref{thm:stefanich-k}]
\label{thm:intro-stefanich}
Higher K-theory defines a functor
\[
 K_{\mathrm{St}}:
 \operatorname{StRing}\longrightarrow
 \operatorname{StRing},
 \qquad
 K_{\mathrm{St}}(\mathbf R)_r
 =K^{[r+1]}(R_{r+1}).
\]
\end{theorem}

The paper is organized as follows.
Section~\ref{sec:conventions} introduces exact sequences
and additive diagrams.
Section~\ref{sec:weights} constructs universal exact sequences.
Section~\ref{sec:presheaf} constructs additive motives by
localizing enriched presheaves, while
Section~\ref{sec:waldhausen} gives an alternative construction
using a weighted Waldhausen $S$-construction on exact flags.
Section~\ref{sec:mapping} establishes the mapping formulas
and comparisons with existing theories.
Section~\ref{sec:higher-k} develops higher K-theory.
Finally, Section~\ref{sec:stefanich} studies the K-theory
of Stefanich rings.

\section*{Acknowledgments}
I am deeply grateful to Seidai Yasuda for his invaluable support during the preparation of this paper.

This work is supported by JSPS KAKENHI Grant Number JP26KJ0476.

\paragraph{Use of AI.}
The author used OpenAI's ChatGPT and Codex as
discussion partners for mathematical arguments,
to assist in revising proofs and their presentation,
and to proofread the manuscript.
The author takes full responsibility for the entire content of this manuscript.

\section*{Notations and assumptions}
\label{sec:notations}
\addcontentsline{toc}{section}{Notations and assumptions}

Every category is an $(\infty,1)$-category. A space is an
$\infty$-groupoid, and $\mathcal S$ denotes the category of small spaces.

Fix nested Grothendieck universes. Small indexing categories belong to
the innermost universe in use. Presentable categories have the
corresponding small colimits, with objects in a larger universe.
Further universes are used for their categories and diagram objects,
without changing the colimit class or compactness of an existing
category.

For a category $C$, $C^\simeq$ is its maximal subgroupoid and
$\Map_C(x,y)$ its mapping space. In a stable category,
$\operatorname{map}_C(x,y)$ is its mapping spectrum.

An object $x$ is compact if $\Map(x,-)$ preserves filtered colimits.
The full subcategory of compact objects is denoted by $A^\omega$.
A presentable category $A$ is compactly generated if $A^\omega$
generates it under filtered colimits.

Write $\mathbb S$ for the sphere, $\Sp$ for the category of spectra, and
$\Sp^\omega$ for the category of finite spectra. For a small idempotent-complete
stable category $C$, $K(C)$ is connective algebraic K-theory, and
$K(R)=K(\Perf(R))$ for a ring spectrum $R$.
Ordinary Ind-completion is $\Ind$. Tensor products of presentable categories are
Lurie's \citep[Section~4.8.1]{LurieHA}.

\section{Exact sequences and additive functors}
\label{sec:conventions}

\subsection{Coefficient modules and their Homs}

\begin{definition}\label{def:presentable-coefficients}
Let $\PrL$ have compactly generated presentable categories as objects
and compact-preserving left adjoints as morphisms. Its tensor product
represents functors preserving colimits separately
\citep[Section~4.8.1]{LurieHA}. We use the symmetric monoidal tower
\begin{equation}\label{eq:coefficient-tower}
 \nPr0=\Sp,\qquad
 \nPr{(r+1)}=\Mod_{\nPr r}(\PrL)\quad(r\geq0).
\end{equation}
Write $\one_r$ for the unit of $\nPr r$.
By \citep[Section~2.1]{Aoki}, the coefficient categories
$\nPr r$ are commutative algebra objects in $\PrL$.
In particular, $\nPr1$ consists of compactly generated
stable presentable categories.
\end{definition}

Ordinary Ind-completion gives the symmetric monoidal equivalence
\begin{equation}\label{eq:stable-ind-equivalence}
 \Ind:\Cat^{\mathrm{perf}}\xrightarrow{\sim}\nPr1,
\end{equation}
where $\Cat^{\mathrm{perf}}$ has small idempotent-complete stable
categories and exact functors. Its inverse takes compact objects
\citep{BGT,LurieHA}.

For $A\in\nPr{(r+1)}$, we write
$\PrHom_A(a,b)\in\nPr r$ for the
coefficient-enriched Hom, characterized by
\begin{equation}\label{eq:tensor-hom}
 \Map_{\nPr r}(v,\PrHom_A(a,b))
 \simeq\Map_A(v\odot a,b).
\end{equation}
Composition and units are supplied by the
coefficient enrichment \citep{GH,Heine}.
At $r=0$ this is the mapping spectrum.

We use the ordinary categories underlying the coefficient modules.
The coefficient-forgetting construction is
\cite[Remark~4.3]{ScholzeGeometry}, with the fixed
$\kappa=\omega$ formalism of
\cite[Definition~1.2 and Section~2.1]{Aoki}.

For $H\in\nPr r$, $r\geq1$,
\[
 \Map_{\nPr r}(\one_r,H)\simeq(H^\omega)^\simeq.
\]
Consequently, for $A\in\nPr{(r+1)}$ and $r\geq1$,
\[
 \Map_A(a,b)\simeq\bigl(\PrHom_A(a,b)^\omega\bigr)^\simeq.
\]
Thus an ordinary arrow in a module is a compact object of its
coefficient Hom. Subsequent local data are read in those coefficient
modules in the same way. We retain their full coefficient Homs.

\subsection{Exact sequences}

\begin{definition}\label{def:adjunction}
Let $H\in\nPr s$, $s\geq2$.
An internal adjunction $i\dashv r$ consists of
arrows $i:a\to b$ and $r:b\to a$, together
with a unit $\eta:1_a\to ri$, a counit
$\varepsilon:ir\to1_b$, and coherent triangle
data in the coefficient Homs.
Equivalently, these are the dual-pair data
of opposite off-diagonal matrices in
\[
 \mathcal B_H(a,b)=
 \prod_{u,v\in\{a,b\}}
 \PrHom_H(u,v)^\omega,
\]
with matrix multiplication given by composition
and finite biproducts.
\end{definition}

Let $H\in\nPr s$, $s\geq2$. Consider adjunctions
$p:b\rightleftarrows c:j$ and $q:b'\rightleftarrows c':k$,
together with arrows $F:b'\to b$ and $G:c'\to c$.
For a $2$-morphism $\beta:pF\Rightarrow Gq$, its right mate is
\[
 Fk\xrightarrow{\eta_p Fk}jpFk
   \xrightarrow{j\beta k}jGqk
   \xrightarrow{jG\varepsilon_q}jG,
\]
where $\eta_p:1_b\to jp$ and $\varepsilon_q:qk\to1_{c'}$
are the unit and counit of the respective adjunctions.
The Beck--Chevalley condition requires this mate to be invertible.

\begin{definition}\label{def:genuine}
For $H\in\nPr1$, an exact sequence is a specified cofiber square
\[
 \begin{tikzcd}
 a\ar[r]\ar[d]&b\ar[d]\\
 0\ar[r]&c.
 \end{tikzcd}
\]
For $H\in\nPr s$, $s\geq2$, an exact sequence
$a\xrightarrow{i}b\xrightarrow{p}c$ consists of $i\dashv r$,
$p\dashv j$, invertible unit $1_a\to ri$, invertible counit
$pj\to1_c$, $pi\simeq0$, and the recursively exact central sequence
\begin{equation}\label{eq:central-row}
 ir\xrightarrow{\varepsilon}1_b\xrightarrow{\eta}jp
       \quad\text{in }\PrHom_H(b,b).
\end{equation}
Its zero comparison is induced by the adjunctions and $pi\simeq0$.
All lower diagram data are retained. We abbreviate the invertible
unit and counit by $ri\simeq1_a$ and $pj\simeq1_c$.

\label{def:exact-sequence-category}
A morphism of complete exact-sequence diagrams has coherent
invertible comparisons for the forward arrows and their
specified right adjoints, compatible with the units and
counits. In particular, the corresponding
Beck--Chevalley mates are required to be invertible.

Write $\EqRow_s(H)$ for this concrete diagram object, including
strong maps of sequences and their coefficient Homs.
The notation $\EqRow_s^\omega(H)$ restricts its three vertices to
$H^\omega$.
\end{definition}

A cofiber square over its first arrow has a contractible space of
choices. With both successive arrows fixed, the specified square
is still part of the data. Adjunction gives $rj\simeq0$ as well.

\subsection{Diagrams lowering the coefficient index}

We write $\Fun^L_{\nPr r}(A,E)$ for the category of
coefficient-linear left adjoints, with strong transformations
and their coefficient Homs. The superscript $L,\omega$
additionally requires preservation of compact objects, while
$\Fun^\oplus$ denotes functors preserving zero objects and
finite biproducts.

Strong transformations have coherent invertible naturality
comparisons. Their coefficient Homs are given by enriched
ends and retain noninvertible intertwiners. The universal
properties below are understood at the level of coefficient
Homs. We do not assume that these functor categories are
compactly generated.

\begin{definition}\label{def:algebraic-diagrams}
For $A\in\nPr{(n+1)}$ and $E\in\nPr n$, define the category $\operatorname{Fun}_n^\omega(A,E)$ by recursion. At $n=1$ it is
$\Fun(A^\omega,E)$. At $n\geq2$ its objects assign $F(a)\in E$
to each $a\in A^\omega$ and local diagrams
\begin{equation}\label{eq:compact-algebraic-diagrams}
 F_{ab}\in\operatorname{Fun}_{n-1}^{\omega,c}
 \bigl(\PrHom_A(a,b),\PrHom_E(Fa,Fb)\bigr),
\end{equation}
with coherent invertible composition and unit comparisons.
Here $\operatorname{Fun}_m^{\omega,c}(B,D)$ is the full subobject
of $\operatorname{Fun}_m^\omega(B,D)$ whose outer values belong to
$D^\omega$. Thus $F(a)$ is arbitrary in $\operatorname{Fun}_n^\omega$,
whereas the local diagrams have compact local values at every index.
Transformation objects and their local Homs are defined by the
corresponding recursive ends.
The notation $\operatorname{Fun}_{n,\mathrm{filt}}(A,E)$ uses all
objects of $A$, with the same local condition, and requires its
ordinary object functor $A\to E$ to preserve filtered colimits.
\end{definition}

Informally, $\operatorname{Fun}_n^\omega(A,E)$ plays the role of
a higher functor category from the underlying $(\infty,n)$-category
of $A$ to that of $E$, with functorial data specified on
compact objects and subject to recursive compact-image
conditions on the coefficient Homs.

\begin{definition}\label{def:additive}
Let $\Add_n^\omega(A,E)$ be the full subobject of
$\operatorname{Fun}_n^\omega(A,E)$ on diagrams preserving zero and
exact sequences. Their images retain the specified
zero comparison and all adjunction data represented at the target
index. Define $\Add_{n,\mathrm{filt}}(A,E)$ analogously.
\end{definition}

We refer to $\Add_n^\omega(A,E)$ as the category of additive
invariants on compact objects, and to
$\Add_{n,\mathrm{filt}}(A,E)$ as the category of
filtered-colimit-preserving additive invariants.
Propositions~\ref{prop:filtered-extension}
and~\ref{prop:global-additive-extension}
show that restriction to compact objects
induces an equivalence between these categories.

Over a rigid symmetric monoidal stable base, these rows are the
split-exact sequences of \cite{BGT}. At the next coefficient level they are the
diagrams used by Mazel-Gee and Stern \citep{BGT,MGS}.
Section~\ref{sec:mapping} gives the normalizations.

\section{Universal exact sequences}
\label{sec:compact-classifiers}\label{sec:weights}

We construct a module classifying exact sequences
directly in the coefficient category.
Ordinary compact objects, duality, and Ind-completion
suffice for its universal property.

\subsection{Coefficient limits and diagram modules}

\begin{lemma}
\label{lem:compact-coefficient-homs}
\label{lem:compact-hom-filtered}
Let $s\geq1$. The forgetful functor
$\nPr s\to\PrL$ creates limits and filtered colimits.
Taking compact objects preserves these operations,
computed in the category of small idempotent-complete
infinity-categories with finite colimits.

For $H=\lim_i H_i$ in $\nPr s$ and $a,b\in H^\omega$,
with images $a_i,b_i$, we have
\[
 \PrHom_H(a,b)
 \simeq\lim_i\PrHom_{H_i}(a_i,b_i).
\]
For a filtered colimit $H=\operatorname*{colim}_i H_i$
and compact $a,b\in H$, choose a common stage $i_0$
representing both objects. Then
\[
 \PrHom_H(a,b)
 \simeq
 \operatorname*{colim}_{(i_0\to j)\in I_{i_0/}}
 \PrHom_{H_j}(a_j,b_j).
\]
Both equivalences are in $\nPr{(s-1)}$.

Moreover, for $H\in\nPr s$ and $a\in H^\omega$,
the functor
$\PrHom_H(a,-):H\to\nPr{(s-1)}$
preserves filtered colimits.
\end{lemma}

\begin{proof}
The first assertion follows from the fixed-cardinal
module formalism \citep[Section~2.1]{Aoki}
and the ordinary Ind-equivalence \citep{LurieHTT}.

For compact $v\in\nPr{(s-1)}$, tensor--Hom
adjunction gives
\begin{equation}\label{eq:fnd-compact-hom-test}
 \Map_{\nPr{(s-1)}}(v,\PrHom_H(a,b))
 \simeq\Map_H(v\odot a,b).
\end{equation}
Since $v\odot a$ is compact, mapping spaces in
the compact cores commute with the indicated
limits and filtered colimits. Compact coefficients
detect equivalences, proving both Hom formulas.

Finally, for a filtered diagram $(b_j)$ in $H$,
compactness of $v\odot a$ gives
\[
 \Map_H(v\odot a,\operatorname*{colim}_j b_j)
 \simeq\operatorname*{colim}_j\Map_H(v\odot a,b_j).
\]
Applying the same tensor--Hom adjunction and
compact detection proves the last assertion.
\end{proof}

A small coefficient category $\mathcal D$ enriched in
$\nPr s$ is an associative algebra of object-indexed coefficient
matrices. Its coherent right modules form
\[
 \Pmod_{\nPr s}(\mathcal D)=\Fun_{\nPr s}(\mathcal D^{\op},\nPr s),
 \qquad h_d(c)=\mathcal D(c,d).
\]
The coefficients in this formula are the presentable objects of
$\nPr s$ themselves. Its parameter-compatible free cocompletion property
is Lemma~\ref{lem:parameter-modules}.

\begin{lemma}\label{lem:parameter-modules}

Let $s\geq1$, let $\mathcal D$ be a small category
enriched in $\nPr{(s-1)}$, and let $E\in\nPr s$.
The diagram module $E^{\mathcal D}$ belongs to
$\nPr s$, with pointwise colimits and coefficient
action. In particular,$\Pmod_{\nPr{(s-1)}}(\mathcal D)
 \simeq \one_s^{\mathcal D^{\op}}\in\nPr s.$
Restriction along the Yoneda embedding gives
\begin{equation}\label{eq:full-presheaf-up}
 \Fun^L_{\nPr{(s-1)}}(\Pmod_{\nPr{(s-1)}}(\mathcal D),E)
 \simeq\Fun_{\nPr{(s-1)}}(\mathcal D,E).
\end{equation}
This equivalence preserves all coefficient intertwining objects.
The represented objects are compact, but their images on the right
may be arbitrary objects of $E$.
\end{lemma}

\begin{proof}
Choose a set of objects $I$ and write $B_{xy}=\mathcal D(x,y)$.
The diagram module $E^{\mathcal D}$ is the category of coherent modules for
\[
 (TF)_y=\coprod_{x\in I}B_{xy}\odot F_x
 \quad\text{on }\prod_{x\in I}E.
\]
Its presentability, pointwise colimits, and free evaluations follow
from \cite[Theorem~4.43 and Corollary~4.44]{HeineBi}. The formulas
\[
 (L_xc)_y=B_{xy}\odot c,\qquad
 \Map(L_xc,F)\simeq\Map_E(c,F_x)
\]
show that compact $c$ gives compact $L_xc$. These objects generate,
and $v\odot L_xc\simeq L_x(v\odot c)$ is compact for compact $v$.
This proves the required compact generation and action condition.
No equality with pointwise compactness of arbitrary diagrams is used.
The coefficient Hom is
\begin{equation}\label{eq:parameter-bar-hom}
 \PrHom_{E^{\mathcal D}}(F,G)
 \simeq\operatorname*{Tot}_{[m]\in\Delta}
 \prod_{y\in I}\PrHom_E((T^mUF)_y,(UG)_y),
\end{equation}
in $\nPr{(s-1)}$, with operators induced by the actions
and the monad
\citep[Corollary~4.49]{HeineBi}. The free cocompletion theorem
\citep[Theorems~3.41(4) and~4.84]{HeineBi} gives
\eqref{eq:full-presheaf-up}, including these coefficient Homs.
\end{proof}

\begin{definition}\label{def:free-linearization}
For $s\geq1$ and a small category $K$, put
\[
 \mathsf F_s(K)=\Fun(K^{\op},\nPr{(s-1)})\in\nPr s.
\]
\end{definition}

Applying Lemma~\ref{lem:parameter-modules} to the free
$\nPr{(s-1)}$-enrichment of $K$ gives a natural equivalence
\begin{equation}\label{eq:free-linearization-full-up}
 \Fun^L_{\nPr{(s-1)}}(\mathsf F_s(K),E)
 \simeq\Fun(K,E),
 \qquad E\in\nPr s,
\end{equation}
including coefficient Homs.
In particular, $\mathsf F_s(*)\simeq\one_s$ and
$\mathsf F_s([1])$ represent compact objects and arrows
between compact objects, respectively. Both modules are
compact in $\nPr s$, since these classification functors
preserve filtered colimits by
Lemma~\ref{lem:compact-coefficient-homs}
and the finiteness of $[1]$.

\begin{remark}\label{rem:parameter-detection}
We detect coefficient Homs by all compact $v\in \nPr{(s-1)}$. Let $J(v)$
have two objects, endomorphism coefficients $\one_{s-1}$, coefficient
$v$ from the first to the second, and zero in the reverse direction.
An object of $E^{J(v)}$ is $(x_0,x_1,v\odot x_0\to x_1)$.
For these finite acyclic parameters, intrinsic compactness is
pointwise. Free evaluations on compact objects have compact values
and generate. Conversely, mapping out of a pointwise compact
diagram is a finite limit of mapping spaces out of compact objects,
using compactness of $v\odot x$, and preserves filtered colimits.
Finite chains detect ordinary functor categories, and these two-object
parameters detect their coefficient intertwining objects.
\end{remark}

\subsection{Exact sequences and their compact classifier}

\begin{lemma}\label{lem:matrix-adjunctions}
For $H\in\nPr s$, $s\geq2$, and an arrow $i:a\to b$ in $H$, coherent
right-adjoint data have an empty or contractible space. Existence
is detected by a reverse arrow, a unit, a counit, and the triangle
identities in the homotopy category of $\mathcal B_H(a,b)$ (see Definition \ref{def:adjunction}).
These spaces preserve limits and filtered colimits with compact
outer endpoints.
\end{lemma}

\begin{proof}
Let $I$ be the matrix in $\mathcal B_H(a,b)$ whose only
nonzero entry is $i$. A right adjoint $r$ of $i$
determines a right dual of $I$, given by the matrix
with $r$ in the opposite off-diagonal entry.
The corresponding unit, counit, and triangle identities
agree with those of the adjunction $i\dashv r$.
Conversely, the diagonal idempotents recover the
opposite off-diagonal entry of any right dual of $I$.
Thus the two constructions give equivalent spaces
of choices.

The existence of a right dual can be checked in the
homotopy category, and its coherent structure is
unique up to a contractible space of choices
\citep[Lemmas~4.6.1.6 and~4.6.1.10]{LurieHA}.
In particular, arbitrary choices of homotopies
witnessing the triangle identities need not themselves
form coherent duality data.

Dual pairs are preserved by limits of monoidal
categories \citep[Proposition~4.6.1.11]{LurieHA}.
For filtered colimits, Lemma~\ref{lem:compact-coefficient-homs}
allows the objects, morphisms, and triangle identities
to be represented at a common stage.
Uniqueness of right adjoints then gives the claimed
compatibility of the spaces of coherent adjunction data.
\end{proof}

\begin{lemma}\label{lem:exact-row-structure}\label{lem:fnd-parameters}

Let $s\geq1$, let $\mathcal D$ be a small category
enriched in $\nPr{(s-1)}$, and let $E\in\nPr s$.
An exact sequence in $E^{\mathcal D}$ consists of
pointwise exact sequences in $E$ whose adjunctions
assemble coherently, satisfying the Beck--Chevalley
conditions along the morphisms of $\mathcal D$.
Morphisms between such exact sequences are pointwise
morphisms compatible with the specified adjunctions,
including their units and counits, and satisfying the
Beck--Chevalley conditions.
Their coefficient Homs are the corresponding
ends. Exactness of a sequence is detected
at its vertices. The functor
\[
 \Phi_s^\omega(H)=\EqRow_s^\omega(H)^\simeq
       \colon\nPr s\longrightarrow\mathcal S
\]
preserves limits and filtered colimits.
\end{lemma}

\begin{proof}
At $s=1$, cofiber sequences are computed pointwise.
For $s\geq2$, pointwise adjunctions assemble precisely
when the Beck--Chevalley mates are invertible
\citep[Proposition~7.3.2.11]{LurieHA}.
The same criterion applies to enriched structure maps,
with coherence inherited from the mate correspondence.
The central exact sequences assemble by induction.
For a map of exact sequences with middle component $f_b:b\to b'$, compare
both central rows in $\PrHom_H(b,b')$ by whiskering. The strong
comparisons identify $f_bir$ with $i'r'f_b$ and $f_bjp$ with
$j'p'f_b$, together with the units and counits. Dual uniqueness
compares their lower adjunction data. 
The specified central zero comparison is canonical by adjunction
and $pi=0$. At the base its cofiber condition is a property of
that specified square. Induction proves the pointwise assertion
and its compatibility with coefficient ends.

In a limit, compact endpoints and their coefficient Homs are
computed by Lemma~\ref{lem:compact-coefficient-homs}. Adjunctions
glue by Lemma~\ref{lem:matrix-adjunctions}, and the central row
glues by the lower-stage assertion. For a filtered colimit,
the same lemmas descend the endpoints, adjunctions, invertible
comparisons, and recursively the central row to a stage. Finite
limits of spaces commute with filtered colimits. Dual uniqueness
and lower induction identify all paths and their compatibilities,
proving the space-valued assertion.
\end{proof}

\begin{proposition}\label{prop:fnd-classifier}\label{def:presentable-exact-classifier}
Let \(s\ge1\).
There is a compact module $W_s\in\nPr s$,
with a distinguished exact sequence of compact objects, such that
\begin{equation}\label{eq:fnd-space-classifier}
 \Map_{\nPr s}(W_s,H)\simeq\Phi_s^\omega(H).
\end{equation}
\end{proposition}

\begin{proof}
The preceding lemma makes $\Phi_s^\omega$ accessible and
limit preserving. The ordinary adjoint functor theorem
\citep[Section~5.5.2]{LurieHTT} gives a left adjoint $F_s$.
Take $W_s=F_s(*)$. Filtered-colimit preservation gives its
compactness. Its distinguished row corresponds to the identity.
\end{proof}

For example, at $s=1$ one may take $W_1=\Fun([1],\Sp)$, with universal exact sequence
at $(\mathbb S\xrightarrow1\mathbb S)$.

\begin{definition}
For \(s\ge 1\), we call \(W_s\) the universal exact-sequence module.
\end{definition}

\begin{corollary}\label{cor:presentable-exact-classifier}
For \(s\ge1\), there are natural equivalences
\begin{align}
 \Fun^L_{\nPr{(s-1)}}(W_s,H)
 &\simeq\EqRow_s(H),\label{eq:presentable-exact-classifier}\\
 \Fun^{L,\omega}_{\nPr{(s-1)}}(W_s,H)
 &\simeq\EqRow_s^\omega(H).\label{eq:compact-exact-classifier}
\end{align}
They preserve all coefficient intertwining objects. The first
equivalence commutes with coefficient diagram parameters.
\end{corollary}

\begin{proof}
First upgrade \eqref{eq:fnd-space-classifier} to the second
assertion. Apply \eqref{eq:fnd-space-classifier} to the finite
parameter modules of Remark~\ref{rem:parameter-detection}, where
compactness is pointwise, and use
Lemma~\ref{lem:fnd-parameters}. Finite chains identify the ordinary
functor categories. The $J(v)$ tests identify every coefficient Hom
by Remark~\ref{rem:parameter-detection}. This proves the complete
compact assertion, without an assumption about arbitrary parameters.

Let us prove the first assertion.
For any small set of objects of $H$, let $C\subset H$
be the full subcategory generated under finite colimits,
retracts, and compact coefficient actions.
Then $K_C=\Ind(C)\in\nPr s$, and the canonical
coefficient-linear left adjoint $K_C\to H$ is
fully faithful on coefficient Homs between objects of $C$.
Indeed, this follows from the tensor--Hom adjunction
and the full faithfulness of $C\subset H$.
Coefficient full faithfulness lifts its entire
recursive data to a compact exact sequence in $K_C$, giving a functor
$W_s\to K_C\to H$. Conversely, an arbitrary coefficient-linear
left adjoint $F:W_s\to H$ lifts through such a $K_C$ after including
$F(W_s^\omega)$ in $C$. Its restriction to $W_s^\omega$ preserves
finite colimits and compact coefficient actions, so ordinary Ind
extends it to a compact-preserving lift. Any small family of
functors, maps of exact sequences, or coefficient intertwiners lifts simultaneously.
Full faithfulness on $C$, tested also after compact coefficient
actions, identifies their mapping spaces and coefficient Homs.
The compact assertion therefore gives the unrestricted assertion.

Let us prove the last assertion.
The equivalences are compatible with coherent families
of diagrams parametrized by small enriched categories,
as modeled by the diagram modules of
Lemma~\ref{lem:parameter-modules}.
In particular, they preserve the coefficient
intertwining objects.
\end{proof}

\subsection{The universal upper exact sequence}

\begin{lemma}\label{lem:canonical-row-maps}
Let \(s\ge1\) and \(H\in \nPr{s}\).
For an exact sequence $e=(a\xrightarrow{i}b\xrightarrow{p}c)$ in \(H\), put
$e_R=(a\xrightarrow1a\to0)$ and
$e_L=(0\to c\xrightarrow1c)$. There is an exact sequence
$e_R\to e\to e_L$ of exact sequences. At $s\geq2$ its
forward maps and their right adjoints have components
\begin{equation}\label{eq:fnd-upper-components}
\begin{array}{c|ccc}
 &a&b&c\\\hline
 e_R\longrightarrow e&1&i&0\\
 e\longrightarrow e_R&1&r&0\\
 e\longrightarrow e_L&0&p&1\\
 e_L\longrightarrow e&0&j&1.
\end{array}
\end{equation}
The construction is natural in complete coefficient parameters.
\end{lemma}

\begin{proof}
At $s=1$ the forward maps give a pointwise cofiber sequence.
At greater $s$, the identities $ri=1$, $pj=1$, $pi=0$, $rj=0$, and the
triangle identities give the displayed strong comparisons and
invertible mates. Whiskering the central row gives
\[
 i\to i\to0,\qquad r\to r\to0,\qquad
 0\to p\to p,\qquad0\to j\to j.
\]
These degenerate rows have contractible compatible completions,
by the cofiber property at the base and dual uniqueness at each
subsequent stage. This establishes their parameter coherences.
The resulting diagram has vertex matrix
\[
 \begin{pmatrix}a&a&0\\a&b&c\\0&c&c\end{pmatrix}.
\]
Every row and column is the given exact sequence or a degenerate sequence.
Lemma~\ref{lem:exact-row-structure} now proves the assertion,
including exactness of the recursive central row.
\end{proof}

The degenerate-sequence functors
$\sigma_L(x)=(x\xrightarrow1x\to0)$ and
$\sigma_R(x)=(0\to x\xrightarrow1x)$ satisfy
\begin{equation}\label{eq:row-evaluation-adjunctions}
 \sigma_L\dashv\operatorname{ev}_a,\qquad
 \operatorname{ev}_c\dashv\sigma_R.
\end{equation}
These adjoints induce equivalences of coefficient Homs.

\begin{proposition}\label{prop:fnd-upper}
Let \(s\ge1\).
There is an internal exact sequence in the coefficient
category $\nPr s$
\begin{equation}\label{eq:presentable-upper-row}
 \one_s\xrightarrow{I}W_s\xrightarrow{L}\one_s.
\end{equation}
Its central exact sequence evaluated at the distinguished middle object is
the universal exact sequence.
\end{proposition}

\begin{proof}
Let $a\to b\to c$ be the universal exact sequence in \(W_s\).
The maps $I,J:\one_s\to W_s$ select $a,c$. The maps
$R,L:W_s\to\one_s$ classify the degenerate sequences respectively $(1_{s-1}\xrightarrow11_{s-1}\to0)$ and
$(0\to 1_{s-1}\xrightarrow11_{s-1})$ in \(1_s\cong \nPr{(s-1)}\). Under
\eqref{eq:presentable-exact-classifier}, precomposition by $I,J$
is vertex evaluation, and precomposition by $R,L$ is
$\sigma_L,\sigma_R$. The coefficient Hom equivalences
\eqref{eq:row-evaluation-adjunctions} give $I\dashv R$ and
$L\dashv J$, with $RI=1$, $LJ=1$, $LI=0$ and $RJ=0$.

Apply Lemma~\ref{lem:canonical-row-maps} to the universal exact sequence.
Transport its exact sequence $e_R\to e\to e_L$ through the
coefficient equivalence \eqref{eq:presentable-exact-classifier}
with $H=W_s$. Its three functors are $IR,1_{W_s},JL$,
and its two maps are the counit and unit just constructed.
Thus $IR\to1_{W_s}\to JL$ is an exact sequence of module
functors. Lemma~\ref{lem:fnd-parameters} makes it pointwise exact.
Evaluation at $b$ gives the universal exact sequence. This uses neither
membership of $\EqRow_s(W_s)$ in a coefficient category nor a
test only at the three distinguished objects.
\end{proof}

For compact $a$ in a coefficient module, $W_s\odot a$ is compact.
Tensor--Hom adjunction represents a compact local exact sequence by a
map $W_s\odot a\to b$. Tensoring
\eqref{eq:presentable-upper-row} with $a$ and evaluating the central
sequence propagates additivity.

\subsection{Ind-completions}
Restriction to compact objects will determine diagrams on the whole
input. Additivity of that extension needs a separate proof.

\begin{proposition}\label{prop:filtered-extension}
Restriction gives a natural equivalence
\[
 \operatorname{Fun}_{n,\mathrm{filt}}(A,E)
 \simeq\operatorname{Fun}_n^\omega(A,E),
\]
compatible with every coefficient-enriched parameter.
\end{proposition}

\begin{proof}
Given $F_C$ on the right, define its object values by the ordinary
filtered Kan extension
\begin{equation}\label{eq:filtered-kan-values}
 F(a)=\operatorname*{colim}_{(c\to a)\in(A^\omega)_{/a}}F_C(c).
\end{equation}
For $n=1$, this is the ordinary Ind universal property.
Assume $n\geq2$. For compact $c$ and
$b=\operatorname{colim}_jb_j$ with compact $b_j$,
Lemma~\ref{lem:compact-hom-filtered} gives
\[
 \PrHom_A(c,b)\simeq\operatorname*{colim}_j\PrHom_A(c,b_j).
\]
Lemma~\ref{lem:compact-coefficient-homs} identifies its compact
objects and their local coefficient data with the filtered colimit
of the corresponding stagewise data. Compose their $F_C$-images
with $F_Cb_j\to Fb$ and use this source colimit to obtain the local
diagram at $(c,b)$.

For $a=\operatorname{colim}_ia_i$ with compact $a_i$, the action
adjunction gives coefficient equivalences
\begin{equation}\label{eq:filtered-first-hom-limits}
 \PrHom_A(a,b)\simeq\lim_i\PrHom_A(a_i,b),\qquad
 \PrHom_E(Fa,Fb)\simeq\lim_i\PrHom_E(F_Ca_i,Fb).
\end{equation}
The compact-core and coefficient-Hom limit statements of
Lemma~\ref{lem:compact-coefficient-homs} apply at every lower index.
Taking the limit of the previously constructed local diagrams defines
the local diagram at $(a,b)$.

For composition, first restrict the source to a compact $a_i$.
Each finite collection of composable local generators and its
comparison maps descends to common compact target stages. At those
stages the asserted comparison is that of $F_C$. Its associativity
and unit identities are the corresponding identities between maps
of associative algebras. One may check them in each simplicial
degree of their free-algebra resolutions. Each degree uses finite
composable strings, so the same descent applies. The coefficient
limits in \eqref{eq:filtered-first-hom-limits} detect the resulting
identities at arbitrary sources.

The construction commutes with coefficient parameters and is determined by
\eqref{eq:filtered-kan-values} and \eqref{eq:filtered-first-hom-limits}.
Apply Remark~\ref{rem:parameter-detection}.
\end{proof}

\begin{proposition}\label{prop:global-additive-extension}
The filtered extension of $F_C\in\Add_n^\omega(A,E)$ belongs to
$\Add_{n,\mathrm{filt}}(A,E)$.
In particular, 
\[
\Add_{n,\mathrm{filt}}(A,E)\cong \Add_n^\omega(A,E).
\]
\end{proposition}

\begin{proof}
Let $F$ be the filtered extension of
$F_C\in\Add_n^\omega(A,E)$.
It preserves zero since $0\in A^\omega$.
For $n=1$, the universal exact sequence
\[
 x\longrightarrow W_1\odot x\longrightarrow x
\]
is a filtered colimit of such sequences with
compact vertices. Its image under $F$ is therefore
exact. The universal property of $W_1$ gives
\[
 F(g)\simeq F(f)+F(h)
\]
for every cofiber sequence $f\to g\to h$
in a coefficient Hom of $A$.
Applying this to $ir\to1_b\to jp$ shows that
the image of any exact sequence
$a\to b\to c$ is split exact.

Now suppose $n\geq2$.
A compact exact sequence in
$\PrHom_A(a,b)$, for compact $a,b$,
is classified by a map $W_n\odot a\to b$.
The universal upper exact sequence
\eqref{eq:presentable-upper-row}, tensored
with $a$, has compact vertices.
Its image under $F_C$ is exact, and
the classifying map recovers the image
of the prescribed local exact sequence.

For arbitrary $a=\operatorname*{colim}_ia_i$
and $b=\operatorname*{colim}_jb_j$ with
compact $a_i,b_j$, use
\[
 \PrHom_A(a_i,b)
 \simeq\operatorname*{colim}_j\PrHom_A(a_i,b_j),
 \qquad
 \PrHom_A(a,b)
 \simeq\lim_i\PrHom_A(a_i,b).
\]
The compact exact-sequence classifier and
Lemma~\ref{lem:compact-coefficient-homs}
show that the local exactness established
above extends to arbitrary $a,b$.

Finally, $F$ preserves the specified adjunctions
and their coherent data. Since it also preserves
the central exact sequence of every input exact
sequence, it preserves exact sequences by
Definition~\ref{def:genuine}.
\end{proof}
\section{Additive motives by presheaf localization}\label{sec:presheaf}

We construct the universal additive invariant by coefficient-linear
presheaf localization.

\subsection{Construction}

\begin{definition}\label{def:raw-presheaves}
For $A\in\nPr{(n+1)}$, define its free presheaf module $\mathsf P_n(A)$ recursively.
Put
\[
 \mathsf P_1(A)=\Fun((A^\omega)^{\op},\Sp).
\]
For $n\geq2$, form the coefficient matrix on $a,b\in A^\omega$
with entries
\begin{equation}\label{eq:raw-coefficient-hom}
 \mathcal R_A(a,b)=\mathsf P_{n-1}(\PrHom_A(a,b))
       \in\nPr{(n-1)},
 \qquad
 \mathsf P_n(A)=\Pmod_{\nPr{(n-1)}}(\mathcal R_A).
\end{equation}
Composition in $\mathcal R_A$ is induced by the lower free external product
and composition in $A$. The external products
\[
 \mathsf P_{n-1}(X)\otimes\mathsf P_{n-1}(Z)
 \longrightarrow\mathsf P_{n-1}(X\otimes Z),
 \qquad X,Z\in\nPr n,
\]
are constructed simultaneously by extending
$y_x\otimes y_z\mapsto y_{x\otimes z}$
for $x\in X^\omega$ and $z\in Z^\omega$
by colimits and coefficient actions.
Associativity and units follow inductively
from the lower free properties.

We also use this construction for a specified small coefficient
matrix $(I,B)$, where $B_{ab}\in\nPr n$ for \(a,b\in I\) and composition and units
are coefficient morphisms. At $n=1$ take spectral presheaves on the
ordinary category with mapping spaces
$\Map_{\nPr1}(\one_1,B_{ab})$.
At $n\geq2$ replace the entries in \eqref{eq:raw-coefficient-hom}
by $\mathsf P_{n-1}(B_{ab})$. Denote the result by
$\mathsf P_n(I,B)$.
\end{definition}

Informally, $\mathsf P_n(A)$ plays the role of
the category of presheaves on the underlying
$(\infty,n)$-category of compact objects of $A$,
with coefficients in $\nPr{(n-1)}$.

Restriction to represented objects gives
\begin{equation}\label{eq:raw-presheaf-up}
 \Fun^L_{\nPr{(n-1)}}(\mathsf P_n(A),E)
 \simeq\operatorname{Fun}_n^\omega(A,E),
 \qquad E\in\nPr n.
\end{equation}
The compact-preserving functors correspond to
$\operatorname{Fun}_n^{\omega,c}(A,E)$.
For $(I,B)$ the same property classifies object values indexed by
$I$ and the corresponding recursive local diagrams.

The free presheaves in Definition~\ref{def:raw-presheaves} implement
this recursion. Each local diagram extends to a compact-preserving
coefficient map
\[
 \mathsf P_{n-1}(\PrHom_A(a,b))
       \longrightarrow\PrHom_E(Fa,Fb).
\]
The compact-image requirement in \eqref{eq:compact-algebraic-diagrams}
is precisely the local compact-preservation condition. It places no
compactness restriction on the outer values $F(a)$.

The presheaf property of Lemma~\ref{lem:parameter-modules} proves
these assertions inductively. Local coefficient maps preserve compact
objects, giving the local compact-image conditions of
Definition~\ref{def:algebraic-diagrams}. Outer object values are
arbitrary in \eqref{eq:raw-presheaf-up}. Applying the construction in
coefficient parameter modules gives the full intertwining objects,
as in Remark~\ref{rem:parameter-detection}.

\begin{definition}\label{def:additive-relations}
Let $A\in\nPr{(n+1)}$. For each complete exact
sequence $e=(a\to b\to c)$ with compact vertices
in $A$, define a map $\theta_e$ in $\mathsf P_n(A)$
recursively.

For $n=1$, put
\[
 \theta_e:
 y_b\amalg_{y_a}y_0\longrightarrow y_c.
\]

For $n=s+1\geq2$, let $e'$ be the central exact
sequence in $\PrHom_A(b,b)$, and write
\[
 \theta_{e'}:x_{e'}\longrightarrow z_{e'}
 \quad\text{in }\mathcal R_A(b,b)
 =\mathsf P_s(\PrHom_A(b,b)).
\]
Put
\[
 J_s=\Fun([1]^{\op},\nPr{(s-1)}),
 \qquad I_s=\one_s,
\]
and let $\varepsilon_s:J_s\to I_s$ be evaluation
at $0$. Since the endpoints of $\theta_{e'}$
are compact, it determines a compact-preserving
coefficient map
\(
 J_s\longrightarrow\mathcal R_A(b,b)\).
By enriched Yoneda and tensor--Hom adjunction,
this corresponds to
\(
 \widetilde\theta_{e'}:
 J_s\odot y_b\longrightarrow y_b\).
Define
\[
 \theta_e:y_b\longrightarrow
 z_e:=y_b\amalg_{J_s\odot y_b}(I_s\odot y_b),
\]
where the other map in the pushout is
$\varepsilon_s\odot y_b$.

For each $m\geq0$, choose a small set
$\mathcal G_m$ representing the compact objects
of $\nPr m$, including the unit. Define
\begin{equation}\label{eq:additive-relation-set}
\begin{split}
 S_n(A)={}&
 \{g\odot(y_0\to0)\mid g\in\mathcal G_{n-1}\}\\
 &{}\cup
 \{g\odot\theta_e\mid
 e\in\EqRow_{n+1}^\omega(A),\
 g\in\mathcal G_{n-1}\}.
\end{split}
\end{equation}
\end{definition}

Choosing representatives of the exact sequences makes the indexing
collections small.

\begin{definition}\label{def:additive-motive}
Define
\begin{equation}\label{eq:presheaf-motive-definition}
 \M_n(A)=\mathsf P_n(A)[S_n(A)^{-1}],
 \qquad L_A:\mathsf P_n(A)\longrightarrow\M_n(A).
\end{equation}
We call $\M_n(A)$ the additive motive of $A$.
\end{definition}

The universality is proved in Theorem \ref{thm:compact-presheaf-up}.

\begin{lemma}\label{lem:compact-module-localization}
Every map in $S_n(A)$ has compact source and target. The localization
\eqref{eq:presheaf-motive-definition} is coefficient-linear, belongs to
$\nPr n$, and preserves compact objects. The objects $g\odot L_Ay_a$,
for $a\in A^\omega$ and compact coefficient $g$, compactly generate its
underlying presentable category.
\end{lemma}

\begin{proof}
The maps $\theta_e$ have compact endpoints
by induction: the base case uses finite colimits
of representables, while the higher cases use
compactness of $J_s$, $I_s$, and their actions
on representables.
The strongly saturated class generated by
$S_n(A)$ is stable under compact coefficient
actions, hence under all coefficient actions
by filtered colimits. Thus the localization
is coefficient-linear.

Since the relations have compact endpoints,
local objects are closed under filtered colimits,
so $L_A$ preserves compact objects.
The objects $g\odot L_Ay_a$, for compact $g$,
are compact generators, proving
$\M_n(A)\in\nPr n$.
\end{proof}

\subsection{The universal property}

\begin{lemma}\label{lem:outgoing-coinverter}
Let $E\in\nPr{(s+1)}$, and let an arrow $\alpha$ of
$\PrHom_E(x,y)$ be classified by $J_s\odot x\to y$. Then
\[
 y\longrightarrow y\amalg_{J_s\odot x}(I_s\odot x)
\]
is an equivalence if and only if $\alpha$ is an equivalence.
\end{lemma}

\begin{proof}
Taking coefficient Hom into $t\in E$ turns the morphism into
\begin{align}\label{equation:evalutaion-map-at-t}
 \PrHom_E(y,t)
 \times_{\PrHom_E(x,t)^{[1]}}\PrHom_E(x,t)\to \PrHom_E(y,t).
\end{align}
The map from the second factor \(\PrHom_E(x,t)\) of the pullback to the arrow category \(\PrHom_E(x,t)^{[1]}\) is the constant-arrow inclusion. Its
image is the equivalence arrows, and each fiber of \(\PrHom_E(x,t)\to \PrHom_E(x,t)^{[1]}\) is either empty or contractible. If $\alpha$ is invertible, every postcomposition lies in that image, so the pullback projection is an equivalence. Testing
against all coefficients and using ordinary Yoneda proves that the
original pushout map is an equivalence. Conversely, assume that the morphism in the statement is an equivalence.
Take $t=y$ and
lift the identity of $y$ along the map in
\eqref{equation:evalutaion-map-at-t}. Its image
in the arrow parameter identifies $\alpha$ with an equivalence arrow.
\end{proof}

\begin{theorem}\label{thm:compact-presheaf-up}\label{thm:presheaf-universal}
Restriction to $a\mapsto L_Ay_a$ gives a natural equivalence
\begin{equation}\label{eq:compact-presheaf-up}
 \Fun^L_{\nPr{(n-1)}}(\M_n(A),E)
 \simeq \Add_n^\omega(A,E)\cong \Add_{n,\mathrm{filt}}(A,E).
\end{equation}
It is compatible with every coefficient-enriched parameter.
\end{theorem}

\begin{proof}
The second equivalence is due to Proposition~\ref{prop:global-additive-extension}.
Let us prove the first equivalence.
By \eqref{eq:raw-presheaf-up}, a coefficient-linear
left adjoint $T:\mathsf P_n(A)\to E$ corresponds to
a diagram $F\in\operatorname{Fun}_n^\omega(A,E)$.
By coefficient-linearity, $T$ inverts the defining
relations \(S_n(A)\) precisely when $F(0)=0$ and $T(\theta_e)$
is an equivalence for every compact exact sequence $e$.

By the universality of localization, it suffices to prove that, under the condition \(F(0_A)=0_E\), \(T(\theta_e)\) is an equivalence if and only if \(F(e)\) is exact.
We prove by induction on $n$.
For $n=1$, this is the cofiber comparison.
For $n=s+1\geq2$, let $e=(a\to b\to c)$ and let
$e'$ be its central exact sequence in $\PrHom_A(b,b)$.
The construction of $\theta_e$ and
Lemma~\ref{lem:outgoing-coinverter} give that \(
 T(\theta_e)\) is an equivalence if and only if \(\Phi_{bb}(\theta_{e'})\) is an equivalence, where
\[
 \Phi_{bb}:\mathsf P_s(\PrHom_A(b,b))
 \longrightarrow\PrHom_E(Fb,Fb)
\]
is the induced coefficient-linear functor.
Its associated diagram preserves zero, since
zero morphisms factor through $0_A$ and $F(0_A)=0_E$.
By induction, the latter condition is equivalent to
exactness of the image of $e'$.
The recursive functoriality of $F$ preserves
adjunctions and their coherent data, so the
image of $e'$ is the central sequence of $F(e)$.
Thus its exactness is equivalent to that of $F(e)$.
\end{proof}

\begin{lemma}\label{lem:module-row-preservation}
The construction $\M_n$ takes internal exact sequences of coefficient
modules and compact-preserving module maps to internal exact sequences
of output modules.
\end{lemma}

\begin{proof}
Parameter functoriality transports the complete duality diagrams,
their units, counits, and mates. On a represented generator $U_H(h)$
of the middle output, the central exact sequence is the image under $U_H$ of
the source central exact sequence evaluated at $h$. It is exact by
Theorem~\ref{thm:presheaf-universal}. The same holds after every
compact coefficient action. All functors involved preserve colimits,
so the base cofiber identities extend from these generators to all
objects. Induction through the transported central diagrams proves
recursive exactness. The proof in every parameter module gives the
same assertion for complete coefficient diagrams.
\end{proof}

\begin{lemma}\label{lem:motive-filtered-colimits}
The functor $\M_n$ preserves filtered colimits
of coefficient modules and compact-preserving
module functors.
\end{lemma}

\begin{proof}
Compact objects and their coefficient Homs
descend along filtered colimits by
Lemma~\ref{lem:compact-coefficient-homs}.
The same holds for compact exact sequences
by Proposition~\ref{prop:fnd-classifier}.
Consequently,
\[
 \Add_n^\omega(\operatorname*{colim}_i A_i,E)
 \simeq\lim_i\Add_n^\omega(A_i,E).
\]
The same argument applies after every coefficient
parameter. Restricting to diagrams with compact
outer values gives the corresponding equivalence
for compact-preserving functors, since compact
objects descend along filtered colimits.
The assertion follows from the universal property
of colimits in the coefficient category.
\end{proof}

\subsection{Tensor products and functoriality}

\begin{proposition}\label{prop:motivic-products}
Compact-preserving coefficient-module maps induce functors on the
$\M_n$, naturally with respect to coefficient parameter diagrams.
There are external products
\[
 \M_n(A)\otimes\M_n(B)\longrightarrow\M_n(A\otimes B)
\]
preserving compact objects, with coherent associativity, symmetry,
and unit. 
In particular, for any $k\geq1$ and any $\mathbb E_k$-algebra
object $A$, Day convolution
endows $\M_n(A)$ with an $\mathbb E_k$-algebra
structure with compact unit $U_A(\one_A)$,
and $U_A$ is $\mathbb E_k$-monoidal.
\end{proposition}

\begin{proof}
The raw external products of
Definition~\ref{def:raw-presheaves}
make $\mathsf P_n$ lax symmetric monoidal.
For compact $b\in B$, the functor
$-\otimes b:A\to A\otimes B$
preserves compact objects and exact sequences.
Thus the composite
\[
 \mathsf P_n(A)\otimes\mathsf P_n(B)
 \longrightarrow\mathsf P_n(A\otimes B)
 \longrightarrow\M_n(A\otimes B)
\]
inverts the defining relations in the first
variable when the second is represented.
The same holds in the other variable.
Since representables and their coefficient actions
generate, the product descends to
\[
 \M_n(A)\otimes\M_n(B)\longrightarrow\M_n(A\otimes B).
\]
It preserves compact generators, and its
coherences follow from those of the raw product
and the parameterized universal property.
The last assertion is a formal fact of lax symmetric monoidal functors.
\end{proof}

\subsection{Comparison with relative BGT motive}

Let $a$ be a small idempotent-complete symmetric
monoidal stable category whose tensor product is
exact in each variable. An $a$-linear stable category
is an object of $\Mod_a(\Cat^{\mathrm{perf}})$.
Its small stable category of linear exact functors
is denoted by
$\Fun_a^{\ex}(\mathcal C,\mathcal D)$.
The base $a$ is called rigid if every object of $a$
is dualizable.

Assume that $a$ is rigid. The category
$\Mod_a(\Cat^{\mathrm{perf}})$ is the category
$\Cat^{\mathrm{perf}}(a)$ of
\citep[Definition~4.3]{HSS}.
Passing to Ind-completions makes it an object of
$\nPr2$, with coefficient Homs
\begin{equation}\label{eq:linear-Hom}
 \PrHom_{\Mod_a(\Cat^{\mathrm{perf}})}
       (\mathcal C,\mathcal D)
 \simeq\Ind\Fun_a^{\ex}(\mathcal C,\mathcal D).
\end{equation}
Indeed, free $a$-modules on compact stable
categories are compact generators, and the
coefficient action preserves compact objects.

\begin{definition}\label{def:relative-additive-motives}
Following \citep[Definition~5.3]{HSS}, a sequence
of $a$-linear stable categories is split-exact
if its underlying sequence is split-exact in
the sense of \citep[Definition~5.18]{BGT}.
Write $\operatorname{Mot}_{\mathrm{add}}^a$
for the category of additive $a$-motives
defined in \citep[Definition~5.8]{HSS}.
Its universal property is given by
\citep[Theorem~5.12]{HSS}.
\end{definition}

\begin{theorem}\label{thm:linear-normalization}
For every rigid $a$ as above, there is a
canonical equivalence
\begin{equation}\label{eq:BGT-normalization}
 \M_1(\Mod_a(\Cat^{\mathrm{perf}}))
 \simeq\operatorname{Mot}_{\mathrm{add}}^a.
\end{equation}
\end{theorem}

\begin{proof}
Since $a$ is rigid, right adjoints of
$a$-linear exact functors are canonically
$a$-linear, by applying
\citep[Proposition~3.11]{BBJ}
after Ind-completion.
By \citep[Proposition~5.4]{HSS}, relative
split-exact sequences are precisely the
localization sequences in the associated
linear $(\infty,2)$-category.
The adjunctions and the central cofiber
sequence $ir\to1\to jp$ identify these
with the exact sequences of
Definition~\ref{def:genuine}.
The equivalence follows from the universal
properties of
Theorem~\ref{thm:presheaf-universal}
and \citep[Theorem~5.12]{HSS}.
\end{proof}

For $a=\Sp^\omega$, this recovers the BGT
additive motive category \citep{BGT}.
The rigidity hypothesis is essential:
without it, an underlying split-exact
sequence need not be internally exact.
\section{The weighted S-construction}\label{sec:waldhausen}

We construct the motives from complete compact flags and fixed
exact-sequence weights. No motive or K-theory functor is applied to
the input Homs.

Fix $n\geq1$ and $\mathcal A\in\nPr{(n+1)}$.

\subsection{Exact flags and their classifying modules}

Let $\operatorname{Ar}[q]$ be the poset of pairs $(i,j)$ with
$0\leq i\leq j\leq q$, ordered coordinatewise.

\begin{definition}\label{def:exact-flags}
For $H\in\nPr d$, with $d\geq1$, write
$\Fun_0(\operatorname{Ar}[q],H)$ for coherent diagrams
$x$ with specified diagonal zeros $x_{ii}=0$.
An \emph{exact flag} of length $q$ is such a
diagram together with a complete exact sequence
\[
 x_{ij}\longrightarrow x_{ik}\longrightarrow x_{jk}
       \qquad(0\leq i<j<k\leq q)
\]
above each three-vertex subdiagram. The sequences share the displayed
objects, arrows, and zero comparisons. Define their structured diagram
data by the pullback
\begin{equation}\label{eq:flag-definition-pullback}
 S_q^{\ex}(H)=
 \Fun_0(\operatorname{Ar}[q],H)
 \mathop{\times}_{\displaystyle\prod_{i<j<k}
                  \Fun_0(\operatorname{Ar}[2],H)}
             \prod_{i<j<k}\EqRow_d(H).
\end{equation}
The maps from $\EqRow_d(H)$ retain the complete coherent square and
forget the extra exactness data. The
formula specifies objects, strong maps, and their local coefficient
Homs separately.

A morphism of flags is a morphism of these complete diagrams. It preserves the complete adjunction diagrams, including the right adjoints and their compatibility with units and counits. Hence the induced Beck–Chevalley mates are equivalences. For two fixed flags
$\xi,\zeta$, let $\mathcal H_{q,H}(\xi,\zeta)\in\nPr{(d-1)}$
be their coefficient Hom. It is the small end of the coefficient
Homs of $H$ over the complete walking diagram, including the
right-adjoint arrows and their unit and counit cells.
\end{definition}

The end is formed in the coefficient category and encodes the full adjunction data, including the right-adjoint comparisons and their Beck--Chevalley compatibilities. 

By Lemma~\ref{lem:exact-row-structure}, it computes precisely the coefficient Homs of strong maps in \eqref{eq:flag-definition-pullback}, rather than merely the maps of the underlying forward diagrams.

Restriction along monotone maps gives $S_\bullet^{\ex}(H)$ a
simplicial structure. Repeated vertices give the canonical degenerate
rows. Compatible right-adjoint data for a fixed left adjoint form a
contractible space by the matrix-duality construction of
Section~\ref{sec:weights}. Hence restriction and composition of flags
give coherent operators. In particular,
$S_0^{\ex}(H)$ is the zero flag and $S_1^{\ex}(H)=H$ as diagram data.

Let $\mathcal F_q(\mathcal A)$ be a small family representing all
complete flags with compact interval objects, together with their
strong maps and coefficient Homs.
For $n\geq2$, form the coefficient matrix
\[
 \mathcal R_{n,q}(\xi,\zeta)
   =\mathsf P_{n-1}
      \bigl(\mathcal H_{q,\mathcal A}(\xi,\zeta)\bigr)
       \in\nPr{(n-1)}.
\]
Composition is induced by the lower free external product and
composition of the ends. Define its raw module by
\begin{equation}\label{eq:waldhausen-input-flags}
 Q_{n,q}(\mathcal A)=
 \begin{cases}
 \Fun(\mathcal F_q(\mathcal A)^{\op},\Sp),&n=1,\\
 \Pmod_{\nPr{(n-1)}}(\mathcal R_{n,q}),&n\geq2.
 \end{cases}
\end{equation}
At $n=1$ the argument of $\Fun$ is the ordinary category of complete
compact flags. At every index the outer represented values are
arbitrary, while local values satisfy the compact-image condition.
The free module property gives its universal property after every
coefficient parameter. In particular $Q_{n,1}(\mathcal A)=
\mathsf P_n(\mathcal A)$. This is the coefficient-matrix variant
of Definition~\ref{def:raw-presheaves}. Restriction gives simplicial
operators preserving represented compact generators.

\begin{definition}\label{def:flag-weights}
Let
\[
 G_{n,q}=\mathsf F_n(\operatorname{Ar}[q])
       \amalg_{\one_n^{\oplus(q+1)}}0,
 \qquad B_n=G_{n,2}.
\]
The map from $\one_n^{\oplus(q+1)}$ selects the diagonal objects.
The universal exact sequence defines $B_n\to W_n$, where $W_n$ is the
module of Proposition~\ref{def:presentable-exact-classifier}. For
$\mathcal T_q=\{(i,j,k)\mid0\leq i<j<k\leq q\}$, set
\begin{equation}\label{eq:flag-weight-pushout}
 T_{n,q}=
 G_{n,q}
 \amalg_{\displaystyle\coprod_{\tau\in\mathcal T_q}B_n}
             \displaystyle\coprod_{\tau\in\mathcal T_q}W_n.
\end{equation}
The other map from each copy of $B_n$ selects the corresponding
three-vertex subdiagram. These are pushouts of coefficient modules.
\end{definition}

\begin{proposition}\label{prop:flag-weight-up}
There are natural parameter-compatible equivalences
\begin{equation}\label{eq:flag-weight-up}
 \Fun^L_{\nPr{(n-1)}}(T_{n,q},E)
       \simeq S_q^{\ex}(E).
\end{equation}
The weights form a cosimplicial module $T_{n,\bullet}$ whose maps
preserve compact objects. Moreover,
\[
 T_{n,0}=0,\qquad T_{n,1}\simeq\one_n,\qquad
 T_{n,2}\simeq W_n.
\]
\end{proposition}

\begin{proof}
The unrestricted exact-sequence classifier
of Corollary~\ref{cor:presentable-exact-classifier}
and the free module property identify
\eqref{eq:flag-weight-pushout} with
\eqref{eq:flag-definition-pullback},
including coefficient parameters.

Restriction of flags induces the cosimplicial
structure. Its maps preserve compact objects
by \eqref{eq:compact-exact-classifier},
since the universal interval objects are compact.
The low-degree identifications follow directly
from the defining pushout.
\end{proof}

At $n=1$, $T_{1,q}$ is the Ind-completion of the perfect category of
a linearly oriented $q$-vertex spectral quiver. Indeed, the usual
gap-diagram equivalence
$S_q^{\ex}(E)\simeq\Fun([q-1],E)$ identifies their representing
properties \citep{Waldhausen}.

\subsection{Recognition from degree one}

The following arguments take place in ordinary presentable categories
and their coefficient Hom objects. They also apply in the parameter
modules of Lemma~\ref{lem:parameter-modules}.

\begin{lemma}\label{lem:exact-row-cofiber}
Let $n\geq 1$ and $E\in\nPr{n}$.
Every internal exact sequence
\(a\xrightarrow{i}b\xrightarrow{p}c\)
in $E$ determines a coherent cofiber square in the underlying presentable category of $E$.
This quotient property is compatible with all strong
parameter diagrams.
\end{lemma}

\begin{proof}
For $n=1$ this is the definition. For $n\geq2$, precomposition by
$p$ is fully faithful, with inverse on its image given by
precomposition by $j$. For $f:b\to z$, postcomposition of the central
exact sequence gives
\[
 fir\longrightarrow f\longrightarrow fjp.
\]
If $fi=0$, the first term is zero.
By induction, the second arrow is an
equivalence.
Thus \(f\) naturally factors through \(p\).
\end{proof}

\begin{lemma}\label{lem:flag-first-row}
The first-row chain
\[
 x_{01}\longrightarrow x_{02}\longrightarrow\cdots
       \longrightarrow x_{0q}
\]
of an exact flag in $E$ determines its interval
quotients and coherent comparisons up to contractible choice. This
holds with every coefficient parameter. A specified map of first-row
chains has either no lift to a morphism of complete flags or a
contractible space of lifts.
\end{lemma}

\begin{proof}
Each $x_{ij}$ is the quotient of $x_{0i}\to x_{0j}$ by
Lemma~\ref{lem:exact-row-cofiber}. The pushout universal property
determines the quotient objects, arrows, and zero comparisons
together. The compatible right-adjoint data are contractible when
they exist.
A map of first-row chains determines a map of the quotient diagrams.
It lifts to a morphism of complete flags exactly when its canonical
right-adjoint mates are equivalences. This is a property of those
comparisons. If it holds, the compatible adjunction data give the
contractible space of lifts. The construction in each parameter
module gives the asserted parameter compatibility.
\end{proof}

For the flag matrices, $\operatorname{Fun}_n
(\mathcal F_p(\mathcal A),S_q^{\ex}(E))$ means assignments of target
flags with coherent unit and composition comparisons and the local
compact-image condition. For $n\geq2$ its local maps use
$\mathsf P_{n-1}(\mathcal H_{p,\mathcal A}(\xi,\zeta))$ and the
target flag coefficient Homs. At $n=1$ these are ordinary functors.
Write
$\operatorname{Fun}_{n,\mathrm{simp}}
(\mathcal F_\bullet(\mathcal A),S_\bullet^{\ex}(E))$
for the end imposing compatibility with simplicial operators,
taken on objects and coefficient Homs.

\begin{theorem}\label{thm:flag-recognition}
Restriction to degree one induces a parameter-compatible equivalence
\begin{equation}\label{eq:flag-recognition}
 \operatorname{Fun}_{n,\mathrm{simp}}
    (\mathcal F_\bullet(\mathcal A),S_\bullet^{\ex}(E))
       \simeq\Add_n^\omega(\mathcal A,E).
\end{equation}
\end{theorem}

\begin{proof}
An additive diagram $F$ preserves exact flags,
including their adjunction data and coherent
comparisons. It therefore induces a simplicial
diagram $S_\bullet^{\ex}(F)$.
Conversely, let $H_\bullet$ be a simplicial
diagram as in the statement, and put $F=H_1$.
Degree zero gives $F(0)=0$.
For an exact sequence $e=(a\to b\to c)$,
the faces identify the vertices of $H_2(e)$
with $F(a),F(b),F(c)$.

By Lemma~\ref{lem:canonical-row-maps}
and \eqref{eq:row-evaluation-adjunctions},
its forward arrows and zero comparison
agree with those of $F(e)$.
Uniqueness of adjunction data
then identifies the remaining structure.
Thus $F$ is additive and
$H_2\simeq S_2^{\ex}(F)$.

For $q\geq3$, the interval faces identify
the first-row chain of $H_q$ with the image
under $F$ of the source first-row chain.
The resulting comparisons are equivalences,
so their right-adjoint mates are invertible.
Lemma~\ref{lem:flag-first-row} therefore
extends them uniquely up to contractible choice
to equivalences $H_q\simeq S_q^{\ex}(F)$.
This uniqueness ensures compatibility with
all simplicial operators.

The same argument in coefficient diagram modules
identifies the coefficient Homs and proves
parameter compatibility.
\end{proof}

\subsection{The weighted realization and its universal property}

\begin{definition}[The weighted Waldhausen \(S\)-construction]\label{def:weighted-s-construction}
With $Q_{n,q}(\mathcal A)$ as in \eqref{eq:waldhausen-input-flags}, define
\begin{equation}\label{eq:weighted-s-construction}
 \operatorname{Wald}_n(\mathcal A)
       =\int^{[q]\in\Delta}T_{n,q}\otimes Q_{n,q}(\mathcal A).
\end{equation}
This coend is the realization of the two-sided simplicial bar
\begin{equation}\label{eq:weighted-s-bar}
 \mathcal B_p=
 \coprod_{[q_0]\to\cdots\to[q_p]\ \mathrm{in}\ \Delta}
       T_{n,q_0}\otimes Q_{n,q_p}(\mathcal A).
\end{equation}
The first face uses the cosimplicial action on $T_{n,\bullet}$.
The last face uses restriction of input flags. The other faces
compose monotone maps, and the degeneracies insert identities.
\end{definition}

The structure maps preserve represented compact generators and
compact coefficient actions, so the tensor products and module
colimit belong to $\nPr n$. Their free properties allow all
coefficient-linear left adjoints in the mapping-out property.
Since $T_{n,1}=\one_n$ and $Q_{n,1}=\mathsf P_n(\mathcal A)$,
the degree-one map gives a
canonical diagram
\[
 U^\omega_{\mathcal A}\in
 \operatorname{Fun}_n^\omega(\mathcal A,
                    \operatorname{Wald}_n(\mathcal A)).
\]

\begin{theorem}\label{thm:weighted-s-universal}
The diagram $U^\omega_{\mathcal A}$ extends to a filtered-continuous additive
diagram on all objects of $\mathcal A$. For every $E\in\nPr n$,
restriction induces a parameter-compatible equivalence
\begin{equation}\label{eq:weighted-s-universal}
 \Fun^L_{\nPr{(n-1)}}(\operatorname{Wald}_n(\mathcal A),E)
       \simeq\Add_{n}^\omega(\mathcal A,E)\cong\Add_{n,\mathrm{filt}}(\mathcal A,E).
\end{equation}
Consequently there is a canonical equivalence
$\operatorname{Wald}_n(\mathcal A)\simeq\M_n(\mathcal A)$, natural
for module morphisms and their coefficient parameter diagrams.
\end{theorem}

\begin{proof}
Let $E^{\mathcal R_{n,p}}$ be the module of diagrams for the raw
flag coefficient matrix, as in Lemma~\ref{lem:parameter-modules}.
At $n=1$ use $\Fun(\mathcal F_p(\mathcal A),E)$.
Free evaluation, tensor currying, and
Proposition~\ref{prop:flag-weight-up} give
\[
\begin{aligned}
 \Fun^L_{\nPr{(n-1)}}
    (T_{n,q}\otimes Q_{n,p}(\mathcal A),E)&\simeq\Fun^L_{\nPr{(n-1)}}(T_{n,q},E^{\mathcal R_{n,p}})\\
 &\simeq S_q^{\ex}(E^{\mathcal R_{n,p}})\\
 &\simeq\operatorname{Fun}_n
       (\mathcal F_p(\mathcal A),S_q^{\ex}(E)).
\end{aligned}
\]
The last step curries complete rows and their right-adjoint data
pointwise in the diagram module.
Mapping out of the coend gives the end
\[
 \Fun^L_{\nPr{(n-1)}}(\operatorname{Wald}_n(\mathcal A),E)
 \simeq\operatorname{Fun}_{n,\mathrm{simp}}
      (\mathcal F_\bullet(\mathcal A),S_\bullet^{\ex}(E)).
\]
Theorem~\ref{thm:flag-recognition} identifies the right side with
$\Add_n^\omega(\mathcal A,E)$.
\end{proof}

\section{Mapping objects and comparison theorems}
\label{sec:mapping}\label{sec:comparison}\label{sec:base}

The categories $\M_n(\mathcal A)$ in this section are those constructed
by localization in Section~\ref{sec:presheaf}. We now compute their
mapping objects. In particular, the formulas below are consequences
of that construction, rather than definitions of the motives.
Ordinary $K$ denotes connective algebraic K-theory of a small
idempotent-complete stable category.

\begin{theorem}\label{thm:mapping-comparison}
Let $n\geq1$ and let $\mathcal A\in\nPr{(n+1)}$.
For $c\in\mathcal A^\omega$ and arbitrary
$a\in\mathcal A$, there are natural equivalences
\begin{align}
 \operatorname{map}_{\M_1(\mathcal A)}(U_{\mathcal A}c,U_{\mathcal A}a)
  &\simeq K(\PrHom_{\mathcal A}(c,a)^\omega)
       \qquad(n=1),\label{eq:mapping-height-one}\\
 \PrHom_{\M_n(\mathcal A)}(U_{\mathcal A}c,U_{\mathcal A}a)
  &\simeq\M_{n-1}(\PrHom_{\mathcal A}(c,a))
       \qquad(n\geq2).\label{eq:mapping-higher}
\end{align}
The first line uses the compact stable category of the coefficient
Hom. In the second line the argument is the presentable coefficient
Hom itself, an object of $\nPr n$.

These equivalences preserve composition and units. They are natural
for compact-preserving module functors and after every coefficient
probe. They also intertwine the external products as lax symmetric
monoidal constructions.
\end{theorem}

We prove the first coefficient level directly. The remaining levels
follow from the already established additive universal properties.

\subsection{The first coefficient level}

Fix $\mathcal A\in\nPr2$ and put $C=\mathcal A^\omega$.
Let $\mathcal X$ be the small stable $2$-category
with objects $C$ and Hom categories
\[
 \mathcal H(c,d)=\PrHom_{\mathcal A}(c,d)^\omega.
\]
Its composition is biexact, and its finite
biproducts are inherited from $\mathcal A$.
The complete exact sequences of
Definition~\ref{def:genuine} agree with the
$1$-localization sequences of
\citep[Definition~1.2.2]{MGS}.
Consequently, the two spectral-presheaf
localizations agree, giving
\[
 \M_1(\mathcal A)
 \simeq\operatorname{Mot}^{\mathrm{MGS}}(\mathcal X).
\]

Write
\[
 \mathscr P_C=\Fun(C^{\op},\Sp),
 \qquad
 y_b=\Sigma^\infty_+\Map_C(-,b),
\]
and let
$L:\mathscr P_C\to\M_1(\mathcal A)$
be the localization.
For $b\in C$, put
\[
 K_b(c)=K(\mathcal H(c,b)).
\]

The proof of \citep[Theorem~2.3.1]{MGS}
shows that $K_b$ is $L$-local and that the
canonical object-to-$K$ map induces an
equivalence
\begin{equation}\label{eq:height-one-unit}
 Ly_b\xrightarrow{\sim}K_b.
\end{equation}
We now establish the coherent multiplicative
form of this comparison.

For a small spectral category $\mathbb B$, write
\[
 \Mod(\mathbb B)=\Pmod_{\Sp}(\mathbb B),
 \qquad
 \Perf(\mathbb B)=\Mod(\mathbb B)^\omega.
\]

\begin{proposition}\label{prop:height-one-presentation}
Let $\mathbb K_C$ be the spectral category
with objects $C$ and mapping spectra
\[
 \mathbb K_C(c,d)=K(\mathcal H(c,d)).
\]
Composition and units are induced by the
lax symmetric monoidal structure of $K$
\citep{BGTMult}.
There is a canonical equivalence
\begin{equation}\label{eq:height-one-presentation}
 \M_1(\mathcal A)\simeq\Mod(\mathbb K_C)
\end{equation}
sending $U_{\mathcal A}b$, for compact $b$,
to the represented module at $b$.
For arbitrary $a\in\mathcal A$, its image
is the module
\[
 c\longmapsto
 K(\PrHom_{\mathcal A}(c,a)^\omega).
\]
The equivalence is natural in $\mathcal A$
and compatible with composition, units,
and external products, including their
coherent finite-arity operations.
\end{proposition}

\begin{proof}
The multiplicative object-to-$K$ map
\citep{BGTMult} induces a spectral functor
\[
 \Sigma^\infty_+C\longrightarrow\mathbb K_C.
\]
Restriction of scalars sends the represented
module $h_b$ to $K_b$. 
By \eqref{eq:height-one-unit}, these
presheaves are local and satisfy
$K_b\simeq Ly_b$.
Since local presheaves
are closed under colimits, restriction
factors through
\[
 R:\Mod(\mathbb K_C)\longrightarrow\M_1(\mathcal A).
\]
On represented modules, the induced map
on mapping spectra is the identity
\[
 \mathbb K_C(c,b)=K_b(c).
\]
Thus $R$ is an equivalence by compact
generation and spectral Yoneda.

For arbitrary $a\in\mathcal A$, the asserted
formula follows by filtered continuity of
$\PrHom_{\mathcal A}(c,-)$ and $K$.
Finally, compatibility with external products
follows from the multiplicativity of the
object-to-$K$ map and colimit preservation,
including the coherent finite-arity operations.
\end{proof}

\subsection{Proof of Theorem~\ref{thm:mapping-comparison}}

\begin{proof}[{Proof of Theorem~\ref{thm:mapping-comparison}}]
We argue by induction on $n$, the case $n=1$
having already been established.
For $n\geq2$, let
\[
 \mathcal D_A^{(n)}(a,b)
 =\M_{n-1}(\PrHom_A(a,b))
\]
and put
\[
 \M'_n(A)=
 \Pmod_{\nPr{(n-1)}}(\mathcal D_A^{(n)}).
\]
Composition and units are induced by the external
products at height $n-1$.
We claim that \(\M_n(A)\simeq\M'_n(A)\).

The natural functor
$A^\omega\to\M'_n(A)$ is additive.
Indeed, applying $\PrHom_A(t,-)$ to an exact
sequence and using
Lemma~\ref{lem:module-row-preservation}
shows that its image is exact after evaluation
at every representable.
The universal property of $\M_n(A)$ therefore
gives a coefficient-linear functor
\[
 \Phi:\M_n(A)\longrightarrow\M'_n(A).
\]

Conversely, the local diagrams of $U_A$ are additive.
Indeed, a compact local exact sequence is classified
by a map $W_n\odot a\to b$, and the universal upper
exact sequence \eqref{eq:presentable-upper-row}
shows that its image under $U_A$ is exact.
By the inductive universal property, these
local diagrams extend to compact-preserving
coefficient maps
\[
 \M_{n-1}(\PrHom_A(a,b))
 \longrightarrow
 \PrHom_{\M_n(A)}(U_Aa,U_Ab).
\]
Their compatibility with composition and units
follows from the parameterized universal property
at height $n-1$. Enriched Yoneda gives
\[
 \Psi:\M'_n(A)\longrightarrow\M_n(A).
\]

Since \(\Psi\Phi U_A=U_A\), \(\Psi\Phi\) is
equivalent to the identity by the universal
property of $\M_n(A)$.
Similarly, $\Phi\Psi$ is the identity on
representables and on their coefficient Homs:
the latter follows by the inductive universal
property from the identity on represented
local generators.
Enriched Yoneda therefore gives
$\Phi\Psi\simeq\mathrm{id}$.

Enriched Yoneda gives the asserted mapping formula
when both objects are compact.
For an arbitrary target
$a=\operatorname*{colim}_i a_i$ with compact $a_i$,
both sides preserve this filtered colimit.
On the left, this follows from the compactness
of $U_Ac$ and the filtered continuity of $U_A$.
On the right, the same follows from
compact-source continuity of
$\PrHom_A(c,-)$ and
Lemma~\ref{lem:motive-filtered-colimits}.
This proves \eqref{eq:mapping-higher}.

Finally, we verify compatibility with external products.
For any finite family $(A_i)_{i=1}^m$, the external
product on the coefficient categories
$\mathcal D_{A_i}^{(n)}$ is induced by
\[
\begin{aligned}
 \bigotimes_i \M_{n-1}(\PrHom_{A_i}(a_i,b_i))
 &\longrightarrow
 \M_{n-1}\Bigl(
   \bigotimes_i\PrHom_{A_i}(a_i,b_i)\Bigr)\\
 &\longrightarrow
 \M_{n-1}\Bigl(
   \PrHom_{\bigotimes_i A_i}
   (\bigotimes_i a_i,\bigotimes_i b_i)\Bigr).
\end{aligned}
\]
By Proposition~\ref{prop:motivic-products},
these maps are coherent in all arities and induce
external products on the enriched presheaf categories
$\M'_n(A)$.

The comparison functors $\Phi$ and $\Psi$ intertwine
these products on representables and their coefficient
Homs, by the inductive universal property.
The parameterized form of that property identifies
all finite-arity operations and their coherences.
Since representables and their coefficient actions
generate under colimits, the compatibility extends
to all objects. Thus the equivalence
$\M_n(A)\simeq\M'_n(A)$ is compatible with
external products as a lax symmetric monoidal
equivalence.

\end{proof}

\section{Higher K-theory}\label{sec:higher-k}

The regular input $\one_{n+1}$ has underlying category $\nPr n$.
Its distinguished compact object is $e_n=\one_n$.

\begin{definition}\label{def:unit-corepresented-K}
Write
$\operatorname{Mot}_n=\M_n(\one_{n+1})$ and
$\xi_n=U_{\one_{n+1}}(e_n)$.
Day convolution has unit $\xi_n$.
For $A\in\nPr n$ define
\begin{equation}
 K^{[n]}(A)=
 \begin{cases}
  \operatorname{map}_{\operatorname{Mot}_1}(\xi_1,U(A)),&n=1,\\
  \PrHom_{\operatorname{Mot}_n}(\xi_n,U(A)),&n\geq2.
 \end{cases}
 \label{eq:unit-corepresented-K}
\end{equation}
\end{definition}

By Proposition~\ref{prop:motivic-products},
the functor $K^{[n]}$ is lax symmetric monoidal.
In particular, for any $\mathbb E_k$-algebra
$A\in\nPr n$, $k\geq1$, the object
$K^{[n]}(A)\in\nPr{(n-1)}$ naturally inherits
an $\mathbb E_k$-algebra structure.

\begin{corollary}\label{cor:K-normalization}
Naturally in $A\in\nPr n$, there are equivalences
\begin{equation}\label{eq:K-normalization}
 K^{[1]}(A)\simeq K(A^\omega)\quad(n=1),
 \qquad
 K^{[n]}(A)\simeq\M_{n-1}(A)\quad(n\geq2).
\end{equation}
These equivalences are lax symmetric monoidal.
In particular, if $A$ is an $\mathbb E_k$-algebra
for $k\geq1$, the equivalences are equivalences of
$\mathbb E_k$-algebras.
\end{corollary}

\begin{proof}
The coefficient Hom from $e_n$ in the regular
module is $A$ itself. Applying
Theorem~\ref{thm:mapping-comparison} gives
the equivalences. Their compatibility with
composition, units, and external products
identifies the lax symmetric monoidal
structures, hence also the induced
$\mathbb E_k$-algebra structures.
\end{proof}

\begin{corollary}\label{cor:compact-internal-hom}
Let $A$ be compact as an object of $\nPr n$ and let $B\in\nPr n$
be arbitrary. Denote by $\operatorname{Func}(A,B)$ the right adjoint
to the tensor product in this specified compact-preserving
coefficient category. Then
\begin{equation}\label{eq:compact-internal-hom}
 \PrHom_{\operatorname{Mot}_n}(U(A),U(B))
 \simeq
 \PrHom_{\operatorname{Mot}_n}
       \bigl(\xi_n,U(\operatorname{Func}(A,B))\bigr).
\end{equation}
For $n=1$ both sides are mapping spectra.
The equivalence is contravariantly natural in compact $A$ and
covariantly natural in $B$.
\end{corollary}

\begin{proof}
In the regular module, tensor--Hom adjunction gives
\[
 \PrHom_{\one_{n+1}}(A,B)\simeq\operatorname{Func}(A,B),
 \qquad
 \operatorname{Func}(e_n,H)\simeq H.
\]
Apply Theorem~\ref{thm:mapping-comparison} first to the compact source
$A$ and then to the compact source $e_n$, with target
$H=\operatorname{Func}(A,B)$.
Both sides of \eqref{eq:compact-internal-hom} identify with
$K^{[n]}(H)$.
\end{proof}

Compactness here is a condition on $A$ in the ambient coefficient
category, not merely compact generation of its underlying category.
For $n=1$ the internal Hom is
$\Ind\Fun^{\ex}(A^\omega,B^\omega)$.
It is not the category of all colimit-preserving functors
$A\to B$.

\begin{corollary}\label{cor:nonzero-all-heights}
For every finite $n\geq1$, the motive $\M_n(\one_{n+1})$ is nonzero.
More precisely,
\[
 \operatorname{map}_{\M_1(\one_2)}(\xi_1,\xi_1)\simeq K(\mathbb S),
 \qquad
 \PrHom_{\M_n(\one_{n+1})}(\xi_n,\xi_n)
       \simeq\M_{n-1}(\one_n)\quad(n\geq2).
\]
\end{corollary}

\begin{corollary}\label{cor:general-MGS}
Let $\mathcal A\in\nPr2$, and let
$\mathcal C_{\mathcal A}$ be the small stable
$2$-category with objects $\mathcal A^\omega$
and Hom categories
\[
 \operatorname{hom}_{\mathcal C_{\mathcal A}}(c,d)
 =\PrHom_{\mathcal A}(c,d)^\omega.
\]
Then there are natural equivalences
\[
 \begin{aligned}
 \operatorname{Mot}^{\mathrm{MGS}}
       (\mathcal C_{\mathcal A})
   &\simeq  K^{[2]}(\mathcal A)
     ,\\
 K^{(2,1)}_{\mathrm{MGS}}
       (\mathcal C_{\mathcal A})
   &\simeq (K^{[2]}(\mathcal A))^\omega,\\
 K^{(2)}_{\mathrm{MGS}}
       (\mathcal C_{\mathcal A})
   &\simeq K((K^{[2]}(\mathcal A))^\omega).
 \end{aligned}
\]
See \cite{MGS} for the notations.
\end{corollary}

\begin{proof}
The first equivalence follows from
Corollary~\ref{cor:K-normalization}
and the comparison with MGS in
Section~\ref{sec:mapping}.
The remaining equivalences follow
from the definitions in \citep{MGS}.
\end{proof}

\section{K-theory of Stefanich rings}
\label{sec:stefanich}

\subsection{Successive K-theories}

For $r\geq0$, the endomorphisms of the monoidal unit
define a functor
\[
 \Omega_r:
 \operatorname{CAlg}(\nPr{(r+1)})\longrightarrow\operatorname{CAlg}(\nPr r),
 \qquad
 C\longmapsto\PrHom_C(\one_C,\one_C).
\]
The commutative algebra structure is induced by the
symmetric monoidal structure of $C$.

\begin{proposition}\label{prop:k-unit-end}
For $n\geq1$ and
$C\in\operatorname{CAlg}(\nPr{(n+1)})$, there is a natural
equivalence of commutative algebras
\begin{equation}\label{eq:k-unit-end}
 \gamma_{n,C}:
 \Omega_{n-1}K^{[n+1]}(C)
 \xrightarrow{\sim}
 K^{[n]}(\Omega_n C).
\end{equation}
\end{proposition}

\begin{proof}
The unit $u=\one_C$ is compact, since the
algebra unit $\one_{n+1}\to C$ preserves
compact objects.
By Corollary~\ref{cor:K-normalization},
the monoidal unit of $K^{[n+1]}(C)$
corresponds to $U_Cu$.
Theorem~\ref{thm:mapping-comparison} therefore
gives
\[
 \begin{aligned}
 \Omega_{n-1}K^{[n+1]}(C)
 &\simeq
 \PrHom_{\M_n(C)}(U_Cu,U_Cu)\\
 &\simeq K^{[n]}(\PrHom_C(u,u))
 =K^{[n]}(\Omega_n C).
 \end{aligned}
\]
For $n=1$, the middle equivalence is the
ordinary $K$-theory mapping formula.
Compatibility with commutative algebra
structures follows from the multiplicative
mapping comparison.
\end{proof}

A \emph{Stefanich ring} in our stable coefficient
tower is a coherent system
\[
 \mathbf R=(R_r,\alpha_r)_{r\geq0},
 \qquad
 R_r\in\operatorname{CAlg}(\nPr r),\qquad
 \alpha_r:R_r\xrightarrow{\sim}\Omega_rR_{r+1}.
\]
Morphisms are coherent families of commutative
algebra maps compatible with the $\alpha_r$.
Write $\operatorname{StRing}$ for the resulting
inverse limit category.
This is the stable, compactly generated variant
of \citep[Definition~3.13]{ScholzeGeometry}.

\begin{theorem}\label{thm:stefanich-k}
Higher K-theory induces a functor
\[
 K_{\mathrm{St}}:
 \operatorname{StRing}\longrightarrow
 \operatorname{StRing},
 \qquad
 K_{\mathrm{St}}(\mathbf R)_r
 =K^{[r+1]}(R_{r+1}).
\]
The structure equivalences are induced by
\eqref{eq:k-unit-end}.
\end{theorem}

\begin{proof}
Corollary~\ref{cor:K-normalization} gives
$K^{[r+1]}(R_{r+1})\in\operatorname{CAlg}(\nPr r)$, and
Proposition~\ref{prop:k-unit-end} gives
coherent equivalences
\[
 \begin{aligned}
 K^{[r+1]}(R_{r+1})
 &\xrightarrow{K^{[r+1]}(\alpha_{r+1})}
 K^{[r+1]}(\Omega_{r+1}R_{r+2})\\
 &\xrightarrow{\gamma_{r+1,R_{r+2}}^{-1}}
 \Omega_r K^{[r+2]}(R_{r+2}).
 \end{aligned}
\]
Naturality yields the asserted functor
on coherent towers.
\end{proof}

In particular, the zeroth level is
\(K_{\mathrm{St}}(\mathbf R)_0
 \simeq K(R_1^\omega)\),
which need not agree with $K(R_0)$.

\subsection{Examples}

By \citep[Lemma~2.11]{Aoki}
(see also \citep[Theorem~3.11]{ScholzeGeometry}),
the functor
\[
 \Mod_{(-)}(\nPr r):
 \operatorname{CAlg}(\nPr r)\longrightarrow
 \operatorname{CAlg}(\nPr{(r+1)})
\]
is fully faithful, with right adjoint $\Omega_r$.
In particular, \(\Omega_r\Mod_A(\nPr r)\simeq A\).

\begin{example}\label{ex:module-stefanich-k}
Let $a$ be a commutative ring spectrum. Define
\[
 R_0(a)=a,\qquad
 R_{r+1}(a)=\Mod_{R_r(a)}(\nPr r).
\]
These form a Stefanich ring, and
\[
 \begin{aligned}
 K_{\mathrm{St}}(\mathbf R(a))_0
   &\simeq K(a),\\
 K_{\mathrm{St}}^2(\mathbf R(a))_0
   &\simeq
   K_{\mathrm{MGS}}^{(2)}
     \bigl(\mathcal C_{R_2(a)}\bigr),
 \end{aligned}
\]
where $\mathcal C_{R_2(a)}$ is the small stable
$2$-category associated with $R_2(a)$ in
Corollary~\ref{cor:general-MGS}.
The second equivalence follows from that corollary
and Corollary~\ref{cor:K-normalization}.

For $a=\mathbb S$, the zeroth level is nonzero,
and hence so are all subsequent levels.
\end{example}

\begin{example}\label{ex:nonmonogenic-stefanich-zero}

Let $k$ be a field of characteristic zero,
viewed as an $\mathbb E_\infty$-ring, and put
\[
 R_0=k,\qquad
 R_1=\Fun(BC_2,\Mod_k),
\]
with pointwise tensor product, where $BC_2$
denotes the classifying space of the cyclic
group of order two.
Then $\Omega_0R_1\simeq k$.
Extend the tower by
$R_{r+1}=\Mod_{R_r}(\nPr r)$ for $r\geq1$.
Since $k[C_2]\simeq k\times k$, we obtain
\[
 K_{\mathrm{St}}(\mathbf R)_0
 \simeq K(k)\oplus K(k),
 \qquad
 K(R_0)\simeq K(k).
\]
Thus the two spectra are not equivalent,
as their $\pi_0$ have different ranks.
\end{example}

\begin{example}
Let $\mathbf R$ be a Stefanich ring.
For every $r\geq0$, the structure equivalence
\[
 K^{[r+1]}(R_{r+1})
 \xrightarrow[\sim]{
   \gamma_{r+1,R_{r+2}}^{-1}
   \circ K^{[r+1]}(\alpha_{r+1})}
 \Omega_rK^{[r+2]}(R_{r+2})
\]
induces, by adjunction, a natural morphism
\[
 \Mod_{K^{[r+1]}(R_{r+1})}(\nPr r)
 \longrightarrow K^{[r+2]}(R_{r+2}).
\]
In particular, the structure equivalence
$\alpha_0:R_0\simeq\Omega_0R_1$ induces
$\Mod_{R_0}(\Sp)\to R_1$, and hence
\[
 K(R_0)\longrightarrow K(R_1^\omega)
 \simeq(K_{\mathrm{St}}(\mathbf R))_0.
\]
Iterating these comparisons gives natural maps
\[
 K^m(R_0)\longrightarrow
 \bigl(K_{\mathrm{St}}^m(\mathbf R)\bigr)_0,
 \qquad m\geq0,
\]
where the superscripts denote iteration.
These maps need not be equivalences.
\end{example}

\subsection{Localizing motives and further directions}

The constructions of this paper suggest a parallel
theory of higher localizing motives. One approach
is to define them recursively, using localizing
motives at the preceding coefficient level.
By Lemma~\ref{lem:outgoing-coinverter},
the resulting
motives admit a presentation as localizations of
$\mathsf P_n(A)$. 
The main remaining issue is to
identify a satisfactory intrinsic notion of higher
Verdier exactness, compatible with the recursive
coefficient structure.

In a recursive Hom-wise model, a localizing analogue of
Theorem~\ref{thm:mapping-comparison} follows formally.
Nisnevich descent is also available. 
These constructions also provide a framework for
higher nonconnective K-theory and its analogues
for coherent Stefanich towers.

Another direction is the construction of higher
traces and Chern characters. In particular, one
may seek higher analogues of Hochschild homology,
topological Hochschild homology, and topological
cyclic homology, together with natural trace and
cyclotomic trace maps from higher K-theory.
A central question is whether these constructions
can retain the relevant noninvertible higher
morphisms and satisfy the required additive or
localizing properties. Such refinements could
provide new computational tools for higher
K-theory and its applications.

\bibliographystyle{amsalpha}
\bibliography{references}
\end{document}